\documentclass[11pt,a4paper]{amsart}
\usepackage{epsf}

\usepackage{amssymb,amsthm,verbatim}

\newtheorem*{notation}{Notation}
\newtheorem{prop}{Proposition}[section]
\newtheorem{thm}[prop]{Theorem}
\newtheorem{lem}[prop]{Lemma}
\newtheorem{cor}[prop]{Corollary}
\newtheorem{conj}[prop]{Conjecture}

\theoremstyle{definition}

\newtheorem{defn}[prop]{Definition}
\newtheorem{rem}[prop]{Remark}

\newtheorem{ques}[prop]{Question}

\newcommand{\ep}{\epsilon}

\newcommand{\Ga}{\Gamma}

\newcommand{\bmat}{\left ( \begin{matrix} }
\newcommand{\emat}{\end{matrix} \right ) }

\newcommand{\ben}{\begin{enumerate}}

\newcommand{\een}{\end{enumerate}}

\newcommand{\blem}{\begin{lem}}

\newcommand{\elem}{\end{lem}}

\newcommand{\bcl}{\begin{clm}}

\newcommand{\ecl}{\end{clm}}

\newcommand{\bthm}{\begin{thm}}

\newcommand{\ethm}{\end{thm}}

\newcommand{\bq}{\begin{ques}}
\newcommand{\eq}{\end{ques}}

\newcommand{\bpr}{\begin{prop}}

\newcommand{\epr}{\end{prop}}

\newcommand{\bco}{\begin{cor}}

\newcommand{\eco}{\end{cor}}

\newcommand{\bcon}{\begin{conj}}

\newcommand{\econ}{\end{conj}}

\newcommand{\bde}{\begin{defn}}

\newcommand{\ede}{\end{defn}}

\newcommand{\bex}{\begin{exa}}

\newcommand{\eexa}{\end{exa}}

\newcommand{\bexe}{\begin{exe}}

\newcommand{\brem}{\begin{rem}}

\newcommand{\erem}{\end{rem}}

\newcommand{\integer}{\mathbb{Z}}

\newcommand{\SR}{\mathsf{SR}_{1 \frac{1}{2}}}
\newcommand{\Gen}{\mathsf{Gen}(\ttt,1)}
\newcommand{\Exp}{\mathsf{Exp}(\ttt,2)}
\newcommand{\ttt}{\mathsf{t}}

\begin{document}

\title{On Kazhdan constants of finite index subgroups in $SL_n(\mathbb{Z})$}

\author{Uzy Hadad}
\date{\today}                   
\maketitle                      

\begin{center}

\end{center}

\noindent \textbf{Abstract.} We prove that for any finite index subgroup $\Ga$ in $SL_n(\mathbb{Z})$, there exists $k=k(n)\in\mathbb{N}$, $\ep=\ep(\Ga)>0$, and
an infinite family of finite index subgroups in $\Ga$ with a Kazhdan constant greater than $\ep$ with respect to a generating set
of order $k$. On the other hand, we prove that for any finite index subgroup $\Ga$ of $SL_n(\mathbb{Z})$, and for any $\ep>0$ and $k\in \mathbb{N}$, there exists a finite index subgroup $\Ga'\leq \Ga$ such that the Kazhdan constant of any finite index subgroup in $\Ga'$ is less than $\ep$, with respect to any generating set of order $k$.
In addition, we prove that the Kazhdan constant of the principal congruence subgroup $\Gamma_n(m)$, with
respect to a generating set consisting of elementary matrices (and their conjugates), is greater than $\frac{c}{m}$, where $c>0$ depends only on $n$.
For a fixed $n$, this bound is asymptotically best possible.

\section{Introduction.}

We begin by recalling the definition of Kazhdan property (T) and Kazhdan constants, first introduced by Kazhdan \cite{Kaz}.

\bde
Let $\Ga$ be a finitely generated discrete group, and let $S \subset \Ga$.
Let $(\pi,\mathcal{H})$ be a unitary representation of $\Ga$ without a non-zero invariant vectors and let $\ep >0$.
A vector $0\neq v \in \mathcal{H}$ is called $(S,\ep)$-invariant, if for any $s \in S$ we have $\|\pi(s)v-v\|\leq \ep \|v\|$.
 Let $\kappa(\Ga,S,\mathcal{H})$ denote the largest number $\ep \geq 0$ such that $\max \limits_{s\in S} \|\pi(s)v-v\|  \geq \ep\|v\|$ for every $v\in \mathcal{H}$. The \textit{Kazhdan constant} $\kappa(\Ga,S)$ is the infimum of the set $\{\kappa(\Ga,S,\mathcal{H})\}$ where $\mathcal{H}$ runs over all the unitary representation without non-zero invariant vector. We say that a group $\Ga$ has \textit{Kazhdan property (T)} if there exists a finite subset $S$ of $\Ga$ such that $\kappa(\Ga,S)>0$.
\ede


In \cite{ShV} Sharma, and Venkataramana showed that any noncocompact irreducible lattice
in a higher real semi-simple Lie group contains a subgroup of finite index which is generated by $3$ elements.
In an analog way, we raised the following question: For a fixed $n\geq 3$, is there exist $\ep >0$ and $k \in \mathbb{N}$,
such that any finite index subgroup  $\Ga$ in $SL_n(\mathbb{Z})$ contains a subgroup of finite index  $\Gamma' \leq \Ga$ with a Kazhdan
constant greater than $\ep$ with respect to some generating set of size $k$. The answer is no, and it follows from the following proposition.

\bpr \label {prop_non_uni_kaz}
Let $\Ga$ be a finite index subgroup in $SL_n(\mathbb{Z})$.
Then for any $\epsilon>0$ and $k\in \mathbb{N}$, there exists a finite index subgroup $\Gamma' \leq \Gamma$, such that
the Kazhdan constant of any finite index subgroup of $\Ga'$ is less than $\ep$, with respect to any set of generators of order $k$.
\epr

On the other hand for a fixed finite index subgroup $\Gamma$ of $SL_n(\mathbb{Z})$ we prove the following:

\bthm \label{thm_uni_kaz} Let $n\geq 3$ and let $\Ga$ be a finite index subgroup of $SL_n(\mathbb{Z})$. Then there exists $\ep=\ep(\Ga)>0$,
 $k=k(n)\in \mathbb{N}$ and an infinite family $\{\Ga_i\}_{i=1}^\infty$ of finite index subgroups of $\Ga$, such that for every $i$, the Kazhdan constant of $\Ga_i$ is greater than $\ep$ with respect to some set of generators of order $k$. 
\ethm

\noindent We remark that in our construction of the family $\{\Ga_i\}_{i=1}^\infty$, the intersection $\bigcap \Ga_i$ is isomorphic to 
some principal congruence subgroup in $SL_{n-1}(\mathbb{Z})$.  \\
\noindent Recall that for any $n\geq 2$, and for any integer $m$, the \textit{principal congruence subgroup} of level $m$ is defined by
$$\Gamma_n(m)=\{g\in SL_n(\mathbb{Z}):g_{ij}-\delta_{ij}\equiv 0 \mbox{ mod } m , 1 \leq i,j\leq n \}.$$

Let $R$ be an associative ring  with a unit and let $i$ and $j$ be distinct integers in
$\mathbb{N}$ with $1 \leq i,j \leq n$. Denote by $X_{ij}(r)$ the
$n \times n$ matrix with $1$ along the diagonal, $r\in R$ in the
$(i,j)$ position, and zero elsewhere. Note that $X_{ij}(-r)$ is the inverse of $X_{ij}(r)$.
 These are the \textit{elementary matrices}. Denote by
$EL_n(R)$ the group which they generate.
Let $R'$ be a subring (not necessarily with a unit) of the ring $R$ . We define
$$EL_n(R')=\langle X_{ij}(r):r \in R', 1 \leq i\neq j \leq n \rangle.$$

\noindent Let $3\leq n \in \mathbb{N}$ and $m \in \mathbb{N}$, set $R_m=m\mathbb{Z}$ and
$S_n(m)=\{X_{ij}(\pm m): 1\leq i\neq j\leq n\}$.
It is obvious that $EL_n(R_m)=\langle S_n(m) \rangle$.
In Section \ref{sec_gen_set_gamma}, we define a generating set $\Sigma_n(m)$ for $\Ga_n(m)$, which contains $S_n(m)$ and a conjugate of $S_n(m)$ by some specific matrix $y_n\in SL_n(\mathbb{Z})$ where $y_n$ depends only on $n$.\\
Shalom \cite{Sh1} and Kassabov \cite{Kas2} showed that the Kazhdan constant $\ep_n$ of $SL_n(\mathbb{Z})$ satisfy $(42\sqrt{n}+860)^{-1}\leq \ep_n < 2n^{-\frac{1}{2}}$ with respect to the elementary matrices $S_n(1)$ as a generating set. Here we study the Kazhdan constant of $\Ga_n(m)$
(respect. $EL_n(R_m)$) with respect to $\Sigma_n(m)$ (respect. $S_n(m)$).

\bthm \label{thm_kaz_el}
Let $n\geq 3$ then there exists a constant $0<c<1$, which depends only on $n$, such that for any $m \in \mathbb{N}$
the Kazhdan constant of $\Ga_n(m)$ $($respect. $EL_n(R_m))$ with respect to the set of generators $\Sigma_n(m)$ $($respect. $S_n(m))$  
is at least $\frac{c}{m}$ and at most $\frac{2}{m}$.
\ethm

\noindent Thus for a fixed $n$, the Kazhdan constant in Theorem \ref{thm_kaz_el} is asymptotically best possible.

\noindent The proofs of the previous theorems use the relative property (T) for the pair $(EL_2(R_m)\ltimes R_l^2,R_l^2)$, and the bounded generation property.

\bde A group $\Ga$ has the \textit{bounded generation} (BG) property if there exist $g_1,...,g_\upsilon \in \Ga$ such that
$$\Ga=\langle g_1\rangle \cdot \cdot \cdot \langle g_\upsilon \rangle.$$
The value $\upsilon$ is called the \textit{degree} of bounded generation.
\ede

Let $X_n(m)$ be the set of all the
elementary matrices in $SL_n(\mathbb{Z})$ that have all off diagonal
entries divisible by $m$. It is clear that $X_n(m) \subset EL_n(R_m)$.
Let $X_n^c(m)=X_n(m)\bigcup X_n(m)^{y_n}$ where $y_n$ is an element in $SL_n(\mathbb{Z})$ which we define in Section \ref{sec_gen_set_gamma}.

The first part in the following theorem is due to D. W. Morris (unpublished).
With his kind permission we have included his argument in section \ref{sec_proof_thm_bg_el}.

\bthm \label{thm_bg_el_n}  Let $n\geq 3$ then there exists a positive integer $\upsilon_n$ such that for all $m\in \mathbb{N}$\\
1) any element in $EL_n(R_m)$ can be written as a product of at most $\upsilon_n$ matrices from $X_n(m)$.\\
2) any element in $\Ga_n(m)$, can be written as a product of at most $2\upsilon_3+n^3$ matrices from $X_n^c(m)$.

\ethm

\noindent We remark that the main difficulty in proving the last theorem is the fact that the ring $R_m$ is a commutative ring without an identity element. Notice also that the degree of the bounded generation ($\upsilon_n$), which is given by the Compactness Theorem, does not yield any explicit bound.\\
In Theorem \ref{thm_kaz_el}, $n$ was fixed and $m$ varied. Now we study the case that $m$ is fixed and $n$ growth to infinity. As a consequence of the the second part of Theorem \ref{thm_bg_el_n} we obtain:

\bco \label{cor_kaz_gamma}
There exists a constant $c>0$, such that for any $m\in \mathbb{N}$ and $n \geq 3 $
$$\kappa(\Ga_n(m),\Sigma_n(m)) > \frac{c}{n^3m}.$$
\eco

\subsection{Kazhdan constants of finitely generated discrete groups}
A basic question in the theory of Kazhdan's property (T) is to compute explicit Kazhdan constants with respect to some generating sets.
This question was raised by Serre and by de la Harp and Valette \cite[p. 133]{HV} and it is has been addressed for several discrete
groups (see the section on Kazhdan constants in the introduction of \cite{BHV}).
The Kazhdan constant is important since it yields quantitative versions for many of the consequences derived from property (T), e.g. expander constant, mixing time on finite graphs and more, we refer to \cite{BHV,Lub} for more details and applications.

In a fascinating paper \cite{Sh1} Shalom (generalized a result of Burger \cite{Bur}) relates property (T) to bounded generation property;
Let $R$ be a commutative, topological unital ring  $R$, and assume that the group $SL_n(R)$ where $n\geq 3$, has
the (BG) property with respect to elementary matrices, then $SL_n(R)$ has Kazhdan property (T) with explicit Kazhdan constant \cite{Sh1}.
Carter and Keller \cite{CK} proved that  $SL_n(\mathcal{O})$ where $\mathcal{O}$ is a ring of integer, has the bounded generation property if $n\geq 3$.

Tavgen \cite{Tav} proved that for $n\geq 3$ the group $EL_n(R_m)$ has the bounded generation property by elementary matrices, but his bound depends on the index $[SL_n(\mathbb{Z}):EL_n(R_m)]$.

\noindent Now if $G$ is a discrete group with Kazhdan constant $\ep=\kappa(G,S)$ where $S$ is a finite generating set of $G$. Then for any finite index subgroup $\Ga \leq G$, there exists a generating set $K$ of $\Ga$ such that $$|K| \leq [G:\Ga]|S|,$$ and $\kappa(\Ga,K) \geq \kappa(G,S)$ (see \cite[Pro. 3.9]{LW}). Since the size of $K$ is increasing with the index of the subgroup, this method is not applicable for the question of estimating Kazhdan constants for an infinite family of finite index subgroups with respect to a bounded set of generators.
Instead, we will use a result by D. W. Morris (unpublished), who proved a bounded generation property for $EL_n(R_m)$ where the bound depends only on $n$.

Another common question in the theory of Kazhdan's property (T), is to find an infinite family of groups
such that any group in this family has Kazhdan constant greater than $\ep$ (for some $\ep>0$) with respect to a bounded set of generators (see the section on Kazhdan constants in the introduction of \cite{BHV}). In this case we say that
these groups has \textit{uniform Kazhdan constants}.
This question is most notable for constructing expander Cayley graphs and it was shown only recently \cite{KLN,BGT} that the family of all finite (non-abelian) simple groups has uniform Kazhdan constants.
In \cite{Had} the author shows that for a fixed ring of integers $\mathcal{O}$, the family $\{SL_n(\mathcal{O})\}_{n=3}^\infty$ has uniform Kazhdan constants.

\subsection{About the proof of Theorem \ref{thm_uni_kaz}} Let $\Ga$ be a finite index subgroup of $SL_n(\mathbb{Z})$ where $n\geq 3$.
From the positive solution to the congruence subgroup problem \cite{Men1,BLS}, $\Ga$ contains $\Ga_n(m)$ for some $m \in \mathbb{N}$.\\
Let $l\in \mathbb{N}$, we show that the relative Kazhdan constant of the pair $(EL_2(R_m)\ltimes R_{m^l}^{2},R_{m^l}^{2})$ depends only on $m$.
Now there are two natural embedding  of $\Ga_{n-1}(m)\ltimes R_{m^l}^{n-1}$ in $\Ga$. One is the case that $R_{m^l}^{n-1}$ is sitting in the lower row and the second is
when $R_{m^l}^{n-1}$ is sitting in the last column. Let $\Ga_n(m,l)$ be the group generated by this two embedding.
We continue by showing that there exists $\delta_n>0$ which depend only on $m$ (and $n$), such that any unitary representation of $\Ga_n(m,l)$ with a $\delta_n$-invariant vector, contains a non-zero $\Ga_n(m^{4l})$-invariant vector. We finish by showing that any element in the quotient group
$\Ga_n(m,l)/\Ga_n(m^{4l})$ can be written as a bounded product of elementary matrices and an element from $\Ga_{n-1}(m)/\Ga_{n-1}(m^{4l})$. Then by the fact the subgroup $\Ga_{n-1}(m)$ has property (T) (or $\tau$ for $n=3$), we get that the Kazhdan constant of $\Ga_n(m,l)$ depends only on $m$ (and $n$).

\subsection{Organization of the Paper.} In Section \ref{pre_kaz} we recall several definition and useful lemmas in the theory of property (T). In Section \ref{sec_prop_non_uni_kaz} we prove Proposition \ref{prop_non_uni_kaz}. In Section \ref{sec_gen_set_gamma} we define a set of generator for $\Gamma_n(m)$, we show that $\Ga_n(m)$ is a bounded product of two subgroups which are isomorphic to $EL_n(R_m)$  and we prove that any element in $EL_n(R_{m^2})$ can be written as a bounded product of elementary matrices from $EL_n(R_m)$ and an element from $\Ga_{n-1}(m)$.
In section \ref{sec_proof_thm_bg_el} we prove Theorem \ref{thm_bg_el_n}.
In section \ref{sec_RT} we prove relative property (T) for the pair $(EL_2(R_m)\ltimes R_l^2,R_l^2)$.
 In Section \ref{sec_kaz_pri_sub} we prove Theorems \ref{thm_kaz_el}.
Finally, in Section \ref{sec_thm_uni_kaz} we prove Theorem \ref{thm_uni_kaz}.

\section{Preliminaries on Kazhdan's property (T)} \label{pre_kaz}

Let $\Ga$ be a discrete group and assume that $\kappa(\Ga,S)>0$ for some finite subset,
then $S$ is a generating set of $\Ga$ (see \cite[Prop. 1.3.2 (ii)]{BHV}). In addition if $\Ga$ has property (T), then $\kappa(\Ga,S)>0$
for any finite generating set $S$ of $\Ga$, but the Kazhdan constant $\kappa(\Ga,S)$ depends on $S$. For the interesting question
of the connection between the Kazhdan constants and the generating set see \cite{ALW,GZ,Os}.

One of the main tools for proving property (T) and estimating Kazhdan constants is the notation of relative property (T).

\bde
Let $\Gamma$ be a discrete group generated by a finite set $S$ and let $B\leq \Ga$ a subgroup. The pair $(\Ga,B)$ is said to have \textit{relative property (T)} if there exists $\ep>0$ such that for every unitary representation $(\pi,\mathcal{H})$ with $(S,\ep)$-invariant vector, contain a non-zero invariant vector under the action of $B$. The largest $\ep$ with this property is called the \textit{relative Kazhdan constant} and is denoted by $\kappa(\Ga,B;S)$.
\ede
It is easy to see that the relative property (T) depends only on the groups $\Ga$ and $B$ and is independent of the generating set $S$ (but the constant depends also on the generating set).\\
\noindent A related useful notation is the Kazhdan ratio given by Ershov and Jaikin \cite{EJ}.
\bde Let $\Gamma$ be a discrete group and $B$ and $S$ subsets of $\Gamma$. The \textit{Kazhdan ratio}
$\kappa_r(G,B;S)$ is the largest $\delta \in \mathbb{R}$ with the following property :
If $\mathcal{H}$ is any unitary representation of $\Gamma$ and $v\in \mathcal{H}$ is $(S,\delta \ep)$-invariant, then $v$ is $(B,\ep)$-invariant.
\ede
\noindent For a discussion about the connection between relative property (T) and Kazhdan ratio see section 3.1 in \cite{EJ}.

If $\kappa(\Ga,B) >0$ for some subset $B$ of $\Ga$ and if $\kappa_r(\Ga,B;S)>0$ for some finite subset
 $S$, then it is easy to check that the group $\Ga$ has Kazhdan property (T) with Kazhdan constant
\begin{equation} \label {eq_kaz_ratio}
\kappa(\Ga,S) \geq \kappa(\Ga,B)\kappa_r(\Ga,B;S).
\end{equation}

\noindent The following (well known) lemmas are useful in proving Kazhdan property (T) and for estimating Kazhdan constants.

\blem \label{lem_kaz_gp} Let $\Gamma$ be a group generated by a set $S$
with $\kappa(\Ga,S) >\ep>0$. Let $0 <\delta <\ep$ be given. Let $(\pi,\mathcal{H})$ be a unitary
representation of $\Gamma$ having a unit vector $v \in
\mathcal{H}$, satisfying $\max \limits_{s \in
S}\|\pi(s)v-v\| \leq \delta.$ Then
$$\max \limits_{g\in \Ga} \|\pi(g)v-v\|<2\frac{\delta}{{\ep}}.$$
 \elem

\noindent For a proof, see for example \cite[Lem. 3.8]{Had}.

\blem \label{cen_mas} Let $(\pi,\mathcal{H})$ be a unitary
representation of a group $\Gamma$. Suppose that for some unit vector $v
\in \mathcal{H}$,  $$\max \limits_{g\in \Gamma} \|\pi(g)v-v\|<\sqrt{2}.$$
 Then there exist a non-zero $\Gamma$-invariant
vector in $\mathcal{H}$.

\elem

\noindent For a proof, we refer the reader to \cite[Ch. \textbf{3}, Cor. \textbf{11}]{HV} and \cite[Lem.\textbf{2.5}]{Sh1}.

The following proposition give a upper bound for the Kazhdan constant of abelian groups.
\bpr \label{pr_kaz_abelian} (\cite[Prop. 3.10]{LW}) Let $A$ be an abelian group of order $n$ generated by $k$ elements $S$, with $\kappa(A,S)=\ep$. Then
$$\ep \leq \frac{2\pi}{|n|^\frac{1}{k}-1}.$$
\epr
\noindent Although the proposition in \cite{LW} state differently (with a different but equivalent definition of Kazhdan constant),
 it follows easily from their proof.

\section{Proof of Proposition \ref{prop_non_uni_kaz}} \label{sec_prop_non_uni_kaz}
Let $\Ga$ be a finite index subgroup of $SL_n(\mathbb{Z})$, and let $\ep \mbox { and } k$ be given. Then
for any prime $p$ the intersection $\Ga'=\Ga_n(p)\bigcap \Ga$ is a subgroup of finite index in $\Ga_n(p)$.
Now the principle congruence subgroup $\Ga_n(p)$ embeds into a
pro-$p$ group, therefore any finite index group in $\Ga_n(p)$, has a quotient isomorphic to the cyclic group $C_p$.
So from Proposition \ref{pr_kaz_abelian} we get that Kazhdan constant of any finite index subgroup in $\Ga'$ is less than $\frac{2\pi}{p^{\frac{1}{k}}-1}$ with respect to any set of generators of order $k$.
Therefore for large $p$, which depends on $\ep$ and $k$, the result follow.
\qed

\section {Generating set for $\Ga_n(m)$ and BG for $\Gamma_n(m^2)$} \label{sec_gen_set_gamma}

\subsection{Generating set}

Let $ 1 \leq i,j\leq n$ and  $E_{ij}$ be the matrix with the entry
$1$ at the $(i,j)$th place and zero elsewhere. For $1 \leq i <n $, let
$$y_n=I-\sum \limits_{i=2}^{n} E_{i+1,i}.$$

e.g. $$y_n=\left(\begin{array}{ccccc}
  1  \\
  -1 & 1 \\
  & & \ddots &  & \\
  &&& 1& 0\\
  &   &  &-1 & 1 \end{array} \right).$$

\noindent Let $A$ be a subset of $SL_n(\mathbb{Z})$ and let $g \in SL_n(\mathbb{Z})$, we define
$A^g=\{a^g:a\in A\}$. Set $S_n(m)=\{X_{ij}(m):1\leq i\neq j \leq n\}$ and $\Sigma_n(m)=S_n(m)\bigcup S_n(m)^{y_n}$.

\bpr \label{pr_bg_princ}Let $m \in\mathbb{N}$ and $n \geq 3$, then the principal
congruence subgroup $\Gamma_n(m)$ in $SL_n(\mathbb{Z})$ is generated by the set $\Sigma_n(m)$ and
$$\Ga_n(m)=EL_n(R_m)\cdot EL_n(R_m)^{y_n}.$$

\epr
\begin{proof}
The group $\Ga_n(m)/\Ga_n(m^2)$ is an abelian group and by direct computation, it is easy to see that any element in $\Ga_n(m)/\Ga_n(m^2)$ is a product of an element from the subgroup $EL_n(R_m)$ and an element from the subgroup $EL_n(R_m)^{y_n}$.
Now, since $\Ga_n(m^2)$ is a normal subgroup of $EL_n(R_m)$ (see \cite{Men2,Tit}) the result follow.

\end{proof}

\subsection{Bounded generation for $\Ga_n(m^2)$} \label{sec_pr_bg_gamma}

We start with an important known property of the ring $\mathbb{Z}$.

\blem \label{lem_sr_z}
Let $a_1,...,a_n$ be non-zero integers which satisfy $\mathbb{Z}a_1+...+\mathbb{Z}a_n=\mathbb{Z}$.
If $n >2$, then there exists $x\in \mathbb{Z}$ such that $$\mathbb{Z}(a_1+xa_n)+\mathbb{Z}a_2+...+\mathbb{Z}a_{n-1}=\mathbb{Z}.$$
\elem
\begin{proof}
We denote by $P$ the set of all primes. Let $d=gcd(a_2,...,a_{n-1})$. Set
$$A=\{p\in P:p|d, p\nmid a_1,p\nmid a_n\}.$$

\noindent If  $A=\emptyset$, set $x=1$  else let $x=\prod \limits_{p \in A} p$, we claims that
$$gcd(a_1+xa_n,a_2,...,a_{n-1})=1.$$

\noindent Indeed, let $p$ be a prime divisor of $gcd(a_2,...,a_{n-1})$. We have three cases:
\begin{enumerate}

\item If $p$ divide $a_1$, then $p \nmid x$ and since $gcd(a_1,...,a_n)=1$
we have $p\nmid a_n$. Hence $$a_1+xa_n\equiv xa_n \mbox{ mod }p.$$

\item If $p$ divide $a_n$ then $p \nmid a_1$ , hence $a_1+xa_n\equiv a_1 \mbox{ mod } p$.

\item If $p\nmid a_1$ and $p\nmid a_n$ then $p\in A$ and $p|x$. Hence $a_1+xa_n\equiv a_1 \mbox{ mod } p$.
\end{enumerate}
\end{proof}

In the language of K-theory, Lemma \ref{lem_sr_z} shows that the stable range of $\mathbb{Z}$ is $2$ (see for example \cite{HO} for the definition of stable range).

\blem \label{lem_sr_ind} Let $n>2$, and let $m,a_1,...,a_n$ be non-zero integers which satisfy $gcd(ma_1,ma_2,...,ma_{n-1},ma_n+1)=1$. 
Then there exists $x_1,...,x_{n-3}\in \mathbb{Z}$ such that
$$gcd(ma_1+x_1m^2a_3+ \cdot\cdot\cdot + x_{n-3}m^2a_{n-1},ma_2,ma_n+1)=1.$$
\elem
\begin{proof}

By Lemma \ref{lem_sr_z}, there exists $x_1\in \mathbb{Z}$ such that
\begin{equation} \label{eq_sr_m}
gcd(ma_1+x_1ma_3,,ma_{4},...,ma_{n-1},ma_2,ma_n+1)=1.
  \end{equation}
\noindent We claim that
  $$gcd(ma_1+x_1m^2a_3,ma_{4},...,ma_{n-1},ma_2,ma_n+1)=1.$$
Indeed, let $d=gcd(ma_{4},...,ma_{n-1},ma_2,ma_n+1)$, and let $p$ be a prime such that $p|d$. It is clear that $p\nmid m$.\\
If $p|x_{1}$, then $p\nmid ma_1$ and we are done.
So assume that $p\nmid x_{1}$ and there are two cases to check.
\begin{enumerate}
\item If $p|ma_{3}$, then from equation \ref{eq_sr_m}, we obtain $p\nmid ma_1$ and we are done.
\item If $p|ma_1$, then from equation \ref{eq_sr_m}, we obtain $p\nmid ma_{i+3}$. This imply that $p\nmid x_{1}m^2a_{3}$ and we are done.
\end{enumerate}

\noindent Now by repeating this argument on the sequence 
$$m(a_1+x_1ma_3),ma_{4},...,ma_{n-1},ma_2,ma_n+1,$$
we obtain the require result.

\end{proof}

\bpr \label{pr_sr_ind}
Let $n>2$ and let $m,a_{1},...,a_n$ be non-zero integers which satisfy $gcd(m^2a_1,...,m^2a_{n-1},m^2a_n+1)=1$.
Then there exists $x_1,...,x_{n-2}\in \mathbb{Z}$ such that
$$gcd(m^2a_{1}+x_1m^3a_3+ \cdot\cdot\cdot + x_{n-3}m^3a_{n-1}+x_{n-2}m(m^2a_n+1),m^2a_{2})=m.$$
\epr

\begin{proof}
From Lemma \ref{lem_sr_ind}, there exists  $x_1,...,x_{n-3}\in \mathbb{Z}$ such that
$$gcd(m^2a_{1}+x_1m^3a_3+ \cdot\cdot\cdot + x_{n-3}m^3a_{n-1},m^2a_{2},m^2a_n+1)=1 .$$
\noindent By Lemma \ref{lem_sr_z} there exists $x_{n-2}\in \mathbb{Z}$ such that
$$gcd(m^2a_{1}+x_1m^3a_1+ \cdot\cdot\cdot + x_{n-3}m^3a_{n-1}+x_{n-2}(m^2a_n+1),m^2a_{2})=1.$$

\noindent and it is clear that $m\nmid x_{n-2}$. Let
$$d=gcd(m^2a_{1}+x_1m^3a_3+ \cdot\cdot\cdot + x_{n-3}m^3a_{n-1}+x_{n-2}m(m^2a_n+1),m^2a_{2}).$$

\noindent It is straightforward to check that $m|d$ but $m^2 \nmid d$, so we can write $d=mx$ for some integer $x$ with $m\nmid x$.
 Assume that there exists a prime $p$ such that $p|x$.\\
If $p|x_{n-2}$ then $p\nmid m^2a_{1}+x_1m^3a_3+ \cdot\cdot\cdot + x_{n-3}m^3a_{n-1}$ and we are done.\\
So assume that $p\nmid x_{n-2}$, and we have two cases:
\begin{enumerate}
\item If $p|m^2a_n+1$ then from the fact that $p|a_2$, we obtain
$$p\nmid m^2a_{1}+x_1m^3a_3+ \cdot\cdot\cdot + x_{n-3}m^3a_{n-1},$$
\noindent  and we are done.
\item If  $$p|m^2a_{1}+x_1m^3a_3+ \cdot\cdot\cdot + x_{n-3}m^3a_{n-1},$$
\noindent then from $p|a_2$ we get $p\nmid  m^2a_n+1$, and therefore $p\nmid x_{n-2}m(m^2a_n+1)$ and the result follow.
\end{enumerate}
\end{proof}

\bpr \label{lem_bg}
Let $n \geq 3$, and $m \in \mathbb{N}$, then any element in $\Gamma_n(m^2)$
can be written as a bounded product of at most $3n-2$  elementary matrices in $EL_n(R_m)$ and an element from $\Gamma_{n-1}(m)$.
\epr

\begin{proof}

Let $$A=\left(\begin{array}{cccc}
  * & \ldots  & *  & m^2a_1\\
    \vdots &  \vdots & \vdots & \vdots \\
  ** & \ldots & * & m^2a_{n-1}  \\
  ** & \ldots & * & m^2a_n+1  \\
\end{array} \right)\in \Gamma_n(m^2).$$

\noindent The matrix $A$ is invertible, therefore $(m^2a_1,...,m^2a_n+1)=1$.
Without restrict the generality we can assume that $a_{1}\neq 0$ and $a_{2} \neq 2$.
Using Lemma \ref{lem_sr_ind}, there exists $x_1,...,x_{n-2}\in \mathbb{Z}$ such that
$$gcd(m^2a_{1}+x_1m^3a_3+ \cdot\cdot\cdot + x_{n-3}m^3a_{n-1}+x_{n-2}m(m^2a_n+1),m^2a_{2})=m.$$

\noindent As a consequence there are $y_1,y_2 \in \mathbb{Z}$ such that
$$y_1(m^2a_{1}+x_1m^3a_1+...+x_{n-2}m(m^2a_n+1))+y_2m^2a_{2}=m,$$
\noindent and therefore
$$y_1m(m^2a_{1}+x_1m^3a_1+...+x_{n-2}m(m^2a_n+1))+y_2m^3a_{2}=m^2,$$
Therefore by performing $n-1$ elementary matrices from $EL_n(R_m)$, i.e.
$$B=X_{n,1}(my_1(-a_n))X_{n,2}(my_2(-a_n))X_{1,n}(x_{n-2}m)\cdot \cdot \cdot X_{1,3}(x_1m)A,$$
we get $1$ in the right lower corner of $B$.
Now by at most $2(n-1)$ elementary matrices we can annihilate all the rest of the last row and the last column.
So by at most $3n-2$ elementary matrices from $X_n(m)$, we can bring $A$ to an element $C$ which belong to $SL_{n-1}(\mathbb{Z})$, but $C= I \mbox{ mod  } m$, hence $C\in \Gamma_{n-1}(m)$.

\end{proof}

\section{Proof of Theorem \ref{thm_bg_el_n}} \label{sec_proof_thm_bg_el}

\noindent \textbf{Proof for part (1) of theorem \ref{thm_bg_el_n}:} This proof is based on a result of Carter, Keller, and Paige \cite{Mo} and is due to D. W. Morris (unpublished).
We start by recalling several notations from \cite{Mo}.
\bde \cite[Def. 2.3]{Mo} Let $X$ be a subset of a group $G$. Let $r$ be a any nonnegative integer, we define $\langle X \rangle_r$, inductively, by:
\begin{itemize}
\item $\langle X\rangle_0=1$,
\item $ \langle X\rangle_{r+1}=\langle X\rangle_r\cdot(X\bigcup X^{-1}\bigcup\{1\})$.
\end{itemize}
 \ede
\begin{notation}
\cite[Def. 2.4]{Mo} Let $A$ be a commutative ring, $\texttt{q}$ be an ideal in $A$, and $n$ a positive integer.
\begin{enumerate}

\item $SL(n,A;\texttt{q})=\{g\in SL(n,A): g \equiv I \mbox{ mod } \texttt{q}\}$.
\item $X_n(\texttt{q})=\{ X_{ij}(a):a\in \texttt{q}, 1\leq i\neq j\leq n \}$.
\item $EL_n(\texttt{q})=\langle X_n(\texttt{q})\rangle$.
\item $X_n^\lhd(A;q)$ is the set of $EL_n(A)$-conjugates of elements of $X_n(q)$.
\item $EL_n^\lhd(A;\texttt{q})=\langle X_n^\lhd(A,\texttt{q})\rangle$. i.e $EL_n^\lhd(A;\texttt{q})$ is the smallest normal subgroup of $EL_n(A)$ that contains $X_n(\texttt{q})$.
\end{enumerate}
\end{notation}

From \cite[Thm.~3.12]{Mo}, we see that if $A$ is any commutative ring satisfying the first-order axioms $\SR$, $\Gen$, and $\Exp$ (see \cite{Mo} for the
definition of these axioms), then
	$$ EL_n^\lhd(A; q^2) \mbox{ is a subgroup of index } \le \ttt \mbox{ in }SL(n,A; q^2)  .$$
Furthermore, \cite[Prop.~2]{Tit} (quoted without proof in \cite[Thm.~6.8]{Mo}) tells us
	$$ EL_n^\lhd(A; q^2)  \subset  \langle X_n(q) \rangle.$$
Therefore
	$$  \langle X_n(q) \rangle \mbox{ contains a subgroup of index }\le \ttt \mbox{ in }SL(n,A; q^2) ,$$
so the Compactness Theorem (cf.\ \cite[Cor.~2.8]{Mo}) implies there exists~$r$ (independent of~$q$) such that
	$$  \langle X_n(q) \rangle_r \mbox{ contains a subgroup of index }\le \ttt \mbox{ in }SL(n,A; q^2) .$$
So
	$$  \langle X_n(q) \rangle_{r+ \ttt} \mbox{ contains }SL(n, A; q^2) \cap  \langle X_n(q) \rangle .$$

\noindent  It is clear that the quotient $ \langle X_n(q) \rangle /SL(n,A;q^2)$
is an abelian group, and since  $X_n(q)$ is a union of $(n-1)^2$ subgroups of elementary matrices, we obtain
$$ \langle X_n(q) \rangle_{r + t+(n-1)^2}  = \langle X_n(q) \rangle.$$

Since $\integer$ satisfies the axioms $\SR$, $\Gen$, and $\Exp$ (see Lemma~2.13, Corollary~3.5, and Theorem~3.9 of \cite{Mo}) with
$\ttt = 2(8!) = 80,640$, we may let $A = \integer$ and the result follows.
\qed\\

\noindent \textbf{Proof for part (2) of theorem \ref{thm_bg_el_n}:} If $n=3$ then from Proposition \ref{pr_bg_princ} we get that any element $\gamma \in\Ga_3(m)$ can be written as $\gamma=\gamma_1\cdot\gamma_2$, where $\gamma_1\in EL_3(R_m)$ and $\gamma_2 \in EL_3(R_m)^{y_3}$. 
Hence by part (1) of this Theorem we get that any element in $\Ga_n(m)$ can be written as a product at most $2\upsilon_3$ matrices from $X_3^c(m)$.

\noindent If $n>3$, set $A_n=\Gamma_n(m)/\Gamma_n(m^2)$. $A_n$ is an abelian group which is generated by the set $\Sigma_n(m)$. It is easy to show that every element in $A_n$
can be written as a product of at most $n^2-1$ matrices from $X_n^c(m)$. The set $A_n$ can be regard also as a set of representatives for the cosets of $\Gamma_n(m)/\Gamma_n(m^2)$. Hence any element $\gamma_n \in \Ga_n(m)$, can be written as $\gamma_n=a_n\gamma_n'$, for some $a_n\in A_n$ and $\gamma_n' \in \Ga_n(m^2)$. By Proposition \ref{lem_bg} we get that $\gamma_n'$ can be written as a product of at most $3n-2$ elementary matrices from $X_n(m)$, and an element from
$\Ga_{n-1}(m)$.

\noindent We repeat this process and get that any element $\gamma_n \in \Gamma_n(m)$
can be written as a product of at most $\sum \limits_{k=4}^n (3k-2)$ elementary matrices in $EL_n(R_m)$ and $n-4$ elements
from the cosets representative $A_n,...,A_4$.
So any element in $\Gamma_n(m)$ can be written by at most
$$\sum \limits_{k=4}^n (3k-2)+\sum \limits_{k=4}^n(k^2-k)\leq n^3$$
matrices from $X_n^c(m)$ and an element from $\Ga_3(m)$.
Now since any element in $\Ga_3(m)$ can we written as a product of at most $2\upsilon_3$ matrices from $X_3^c(m)$,
we obtain that any element in $\Ga_n(m)$ can be written as a product at most $2\upsilon_3+n^3$ matrices from $X_n^c(m)$.
\qed

\section{Relative property (T) for $(EL_2(R_m)\ltimes (R_l)^2,(R_l)^2)$} \label{sec_RT}

Burger \cite{Bur} gave the first quantitative version of relative property (T) for the pair of groups $(SL_2(\mathbb{Z}\ltimes \mathbb{Z}^2,\mathbb{Z}^2 )$,
Shalom \cite{Sh1} extended it to the pair $(EL_2(R)\ltimes R^2,R^2)$ where $R$ is an commutative ring with a unit, and Kassabov \cite{Kas1}
proved an analog version for the case that $R$ is an associative ring with unit. Here we give a similar result for some commutative rings
without a unit element. Our proof is based on Shalom proof \cite[Thm. 12.1]{Sh1} which appear also in \cite[Ch. 4]{BHV}.

Let $m,l \in \mathbb{N}$, set $$U^{\pm}(m)=\left(\begin{array}{cc}
  1 & \pm m  \\
  0 & 1\end{array} \right)
  \mbox { and  } L^{\pm}(m)=\left(\begin{array}{cc}
  1 & 0  \\
  \pm m & 1\end{array} \right),$$
in $EL_2(R_m)$, and let  $$e^{\pm}(l)=\left(\begin{array}{c}
  \pm l  \\
  0 \end{array} \right)
  \mbox { and  } f^{\pm}(l)=\left(\begin{array}{c}
  0  \\
  \pm l\end{array} \right),$$ in $R_l^2=({l\mathbb{Z}})^2$.

It is easy to check that the
  set $$F(m,l)=\{U^{\pm},L^{\pm},e^{\pm},f^{\pm}\}$$
  is a generating set for the semidirect product
  $EL_2(R_m) \ltimes (R_l)^2$.

\blem \label{Lem_inv_measure} Let $\nu$ be a mean on the Borel sets of $\mathbb{R}^2 \setminus \{0\}$. There exists a Borel subset $M$ of $\mathbb{R}^2\setminus \{0\}$ and an elementary element $\gamma\in \{U^{\pm}(m),L^{\pm}(m)\}$ such that $|\nu(\gamma(M)-\nu(M)|\geq \frac{1}{4}$ for the linear action of $EL_2(R_m)$ on $\mathbb{R}^2$. \elem

In \cite[Lem 4.2.1]{BHV} the above lemma is stated for the case $m=1$, but it easy to check that the same proof holds for arbitrary $m\in \mathbb{N}$.

The abelian group $R_l^2$ is isomorphic to $\mathbb{Z}^2$ and therefore the unitary dual of $R_l^2$
will be identified with the $2$-torus $\textbf{T}^2$ by associating to $(e^{2\pi ix},e^{2\pi iy})\in \textbf{T}^2$ the character
$(al,bl)\rightarrow e^{2\pi i(xa+yb)}$. The unitary dual $\textbf{T}^2$ is a locally compact space and we denote by $\beta(\textbf{T}^2)$ its Borel subsets. Set $\ep(m)=\frac{\sin\frac{\pi}{m}}{27}$.

The next theorem is quite standard in the theory of Kazhdan's property (T). For a proof see for example Theorem 4.2.2 in \cite{BHV} (just replace $\frac{1}{10}$ with $\ep(m)$ and use $X=(-\frac{1}{2(m+1)},\frac{1}{2(m+1)}]^2$).

\bthm \label{thm_rt_el} Let $m,l \in \mathbb{N}$, then the pair $(EL_2(R_m)\ltimes R_l^2,R_l^2)$ has relative property (T) with relative Kazhdan constant $$\kappa(EL_2(R_m)\ltimes R_l^2,R_l^2;F(m,l))\geq \ep(m).$$  \ethm

The following proposition is analog to Corollary \textbf{4.2.3} in \cite{BHV}. We refer the reader to their proof (just to replace $\mathbb{Z}$ by $R_m$ and $\frac{\ep}{20}$ by $\frac{\delta\ep(m)}{2}$, all the rest is the same).

\bpr \label {rt_SL_2} Let $m,l \in \mathbb{N}$, then  $$\kappa_r(EL_2(R_m)\ltimes R_l^2,R_l^2;F(m,l))>\frac{\ep(m)}{2}.$$
\epr

In case $m=l$ we can use the natural embedding of $EL_2(R_m)\ltimes R_m^2$ into $EL_n(R_m)$ where $n\geq 3$ (see for example \cite[Lem. 3.2.4]{BHV}), and conclude the following proposition (which is analog to \cite[Cor. 1.10]{Kas1}):

\bpr \label{pro_kaz_ratio_el} Let $n\geq 3$ and $m\in \mathbb{N}$, then
$$\kappa_r(EL_n(R_m),X_n(m);S(m)) > \frac{\ep(m)}{2}.$$
\epr

\section{Proof of Theorem \ref{thm_kaz_el} and Corollary \ref{cor_kaz_gamma}} \label{sec_kaz_pri_sub}

\noindent \textbf{Proof of Theorem \ref{thm_kaz_el}:} We start by proving the lower bound separately for each group.\\
\noindent \textbf{The group $EL_n(R_m)$:}
From Theorem \ref{thm_bg_el_n} and Lemma \ref{cen_mas} we get $$\kappa(EL_n(R_m),X_n(m))>\frac{1}{\upsilon_n}.$$
So by Corollary \ref{pro_kaz_ratio_el} and inequality (\ref{eq_kaz_ratio}) we obtain
 $$\kappa(EL_n(R_m),S_n(m)) >  \frac{1}{\upsilon(n)}\cdot\frac{\ep(m)}{2}=\frac{\ep(m)}{2\upsilon(n)}.$$

\noindent \textbf{The group $\Ga_n(m)$:} From Proposition \ref{pr_bg_princ} we get that any element in the group $\Ga_n(m)$ can be written as a  product of $2$
elements from $2$ subgroup which are isomorphic
to $EL_n(R_m)$. Now by the previous case we have $$\kappa(EL_n(R_m),S_n(m))> \frac{\ep(m)}{2\upsilon(n)}.$$
Hence by Lemma \ref{lem_kaz_gp} (and Lemma \ref{cen_mas}) and from the structure of $\Sigma_n(m)$ we obtain
$$\kappa(\Ga_n(m),\Sigma_n(m))\geq  \frac{1}{4}\cdot \frac{\ep(m)}{2\upsilon(n)}=\frac{\ep(m)}{8\upsilon(n)}.$$

\noindent \textbf{The upper bound:} It easy to check that $[\Ga_n(m),\Ga_n(m)] \subset \Gamma_n(m^2)$, and it is known that $\Gamma_n(m^2) \leq EL_n(R_m)$ (see \cite{Tit,Men2}). This imply that the quotient $\Ga_n(m)/\Ga_n(m^2)$ (respect. $EL_n(R_m)/\Gamma_n(m^2)$) is an abelian group which is isomorphic to $C_m^{n^2-1}$
(respect $C_m^{n^2-n}$) where
$C_m$ is the cyclic group of order $m$. From the structure of $\Sigma_n(m)$ (respect $S_n(m)$), one can show that there exists a mapping from $\Ga_n(m)$
(respect $EL_n(R_m)$) onto the cyclic group $C_m$ where the projection of the generating set $\Sigma_n(m)$ (respect. $S_n(m)$) is $\{\pm 1\}$. Now it is straightforward to check that the Kazhdan constant of the cyclic group $C_m$ with respect to $\{\pm1\}$ as generating set is $2\sin\frac{\pi}{m}\sim \frac{2}{m}$.
\qed\\

\noindent \textbf{Proof of Corollary \ref{cor_kaz_gamma}:} 
From part (2) of Theorem \ref{thm_bg_el_n} and Lemma \ref{cen_mas}, we obtains $\kappa(\Ga_n(m),X_n^c(m))>\frac{1}{2\upsilon_3+n^3}$.
By the structure of the generating set $\Sigma_n(m)$ (and Proposition \ref{pro_kaz_ratio_el}), it is straightforward to check

$$\kappa_r(\Ga_n(m),X_n^c(m);\Sigma_n(m))\geq \frac{\ep(m)}{2}.$$

\noindent Hence by equation \ref{eq_kaz_ratio} the result follow.

\qed

\section{Proof of Theorem \ref{thm_uni_kaz}} \label{sec_thm_uni_kaz}

Let $m,l \in \mathbb{N}$ and $n \geq 3$, set

$$\alpha_n(m,l)=\{X_{in}(m^l):1 \leq i < n\}\bigcup \{X_{nj}(m^l):1\leq j < n\}.$$

\noindent If $n >3$ set

$$\beta_n(m,l)=\Sigma_{n-1}(m)\bigcup \alpha_n(m,l),$$
where $\Sigma_{n-1}(m)$ is the set of generators of $\Ga_{n-1}(m)$ which we defined in section \ref{sec_gen_set_gamma}.\\
If $n=3$ let $\Sigma_2(m)$ be any finite set of generators for $\Ga_2(m)$ which contains $\{X_{12}(m),X_{21}(m)\}$ ($\Ga_2(m)$
is of finite index in $SL_2(\mathbb{Z})$, so such a generating set exists), and set
 $$\beta_3(m,l)= \alpha_3(m,l)  \bigcup \Sigma_2(m).$$

Let $\Ga_n(m,l) $ be the group which is
generated by $\beta_n(m,l)$. From a result of Tits \cite{Men2,Tit}, it is easy to check that  $\Ga_n(m,l)$ is a finite index subgroup of  $SL_n(\mathbb{Z})$.

\noindent Set
$$X_n'(m^{l})=\{X_{ij}(t): t\in R_{m^{l}},3 \leq i\leq n, 1\leq j <n\}\bigcup \{X_{lk}(t):t\in R_{m^{l}},3\leq k\leq n,1\leq l <n\}.$$

\blem \label{lem_el_2_emb} Let $g\in X_n'(m^{l})$, then there exists an injective homomorphisms
 $$\psi:EL_2(R_m)\ltimes R_{m^l}^2 \rightarrow \Ga_n(m,l)$$ such that $g$ is contained in $\psi(R_{m^l}^2)$.

\elem
\noindent For a proof see Lemma \cite[Lem. 4.2.4]{BHV}.
The next proposition follow from the last Lemma and Proposition \ref{pro_kaz_ratio_el}.

\bpr \label{pro_kaz_ratio_ga} Let $m,l \in \mathbb{N}$, then
$$\kappa_r(\Ga_n(m,l),X_n'(m^{l}) ;\beta_n(m,l))\geq \frac{\ep(m)}{2}.$$
\epr

Let $\upsilon_n=\upsilon(n)$ be the degree of $EL_n(R_{m^{2l}})$  (recall that $\upsilon_n$ depends only on $n$), and let $\delta_n=\frac{\ep(m)}{8\upsilon_n}$. In the next proposition we use the following obvious observation  $$EL_n(R_{m^{2l}})\leq \Ga_n(m,l).$$

\bpr  \label{prop_inv_el_m_l} Let $(\pi,\mathcal{H})$ be a unitary representation of $\Ga_n(m,l)$. Assume that there exists a unit vector $0\neq v \in \mathcal{H}$ which is $(\beta_n(m,l),\delta_n)$-invariant. Then there exists a non-zero vector in
$\mathcal{H}$ which is invariant under the action of $EL_n(m^{2l})$.
\epr

\begin{proof}

From the Steinberg commutator relations we have $[X_{in}(m^l),X_{nj}(m^l)]=X_{ij}(m^{2l})$. Hence any element in $X_n(m^{2l})$ can be written as
a product of at most $4$ elements from $X_n'(m^l)$, and by the bounded generation property
for $EL_n(R_{m^{2l}})$ we get that any element in $EL_n(R_{m^{2l}})$ can be written as a product of at most $4\upsilon_n$ elementary matrices from $X_n'(m^{l})$.

\noindent Now let $(\pi,\mathcal{H})$ be a unitary representation of $\Ga_n(m,l)$. Assume that there exists a unit vector $0\neq v \in \mathcal{H}$ which is $(\beta_n(m,l),\delta_n)$-invariant. From the above we get that for any $g \in EL_n(m^{2l})$ there exits $g_1,...,g_{4\upsilon_n}\in X_n'(m^l)$
such that
$$g=g_1\cdot \cdot\cdot g_{4\upsilon_n}.$$

\noindent Hence from Proposition \ref{pro_kaz_ratio_ga} we obtain that
$$\|\pi(g)v-v\| \leq \sum\limits_{i=1}^{4\upsilon_n}\|\pi(g_i)v-v\|\leq 4\upsilon_n\frac{\delta_n}{\frac{\ep(m)}{2}}\leq 4\upsilon_n\frac{1}{4\upsilon_n}\leq 1,$$
and from Lemma \ref{cen_mas} the result follow.

\end{proof}

\noindent Let $G_n(m,l)=\Ga_n(m,l)/\Ga_n(m^{4l})$ and let $K$ be the kernel of the following natural homomorphism
$$K=\ker (G_n(m,l)\rightarrow SL_n(\mathbb{Z}/m^{2l}\mathbb{Z})).$$

\blem Any element in $G_n(m,l)/K$ can be written as a product  of at most $2(n-1)$ elementary matrices from $X_n'(m^l) (\mbox{ mod } m^{2l})$ and an
element from $\Ga_{n-1}(m) (\mbox{ mod }m^{2l})$.
\elem

\begin{proof}
Let $A$ be the subgroup generated by $\alpha_n(m,l) (\mbox{ mod } m^{2l})$. It is easy to check that $A$ is an abelian normal subgroup of $G_n(m,l)/K$ and any element in $A$ can be written as a product of at most $2(n-1)$ elementary matrices from $X_n'(m^l) \mbox{ mod }  m^{2l}$.
Moreover, the quotient $(G_n(m,l)/K)/A$ is isomorphic to $\Ga_{n-1}(m)/\Ga_{n-1}(m^{2l})$. Hence any element in $G_n(m,l)/K$ can be written as a product of at most $2(n-1)$ elementary matrices from $X_n'(m^l) (\mbox{ mod }  m^{2l})$ and an element from $\Ga_{n-1}(m) (\mbox{ mod } m^{2l}) $.
\end{proof}

\blem Any element in $K$ can be written by at most $3n-2$ matrices from $X_n'(m^l)$ $\mbox{ mod } m^{4l}$ and an element from  $\Ga_{n-1}(m)/\Ga_{n-1}(m^{4l})$.
\elem

\begin{proof}
It is clear that $K$ is a subgroup of $\Ga_n(m^{2l})/\Ga_n(m^{4l})$.  Let $g \in K$, by Proposition \ref{lem_bg} we obtain that
$g$ can be written by at most $3n-2$ matrices from $X_n'(m^l) (\mbox{ mod } m^{4l})$ and an element from  $\Ga_{n-1}(m)( \mbox{ mod } m^{4l})$ (which is
isomorphic to $\Ga_{n-1}(m)/\Ga_{n-1}(m^{4l})$).

\end{proof}

\noindent From the last two lemmas it is straightforward to check the following:

\bpr \label{ga_mod_bg}
Let $g$ be an arbitrary element in $G_n(m,l)$, then there exists $x_1,...,x_{5n-4} \in X_n'(m^l)(\mbox{ mod } m^{4l})$ and $\gamma_1,\gamma_2\in \Ga_{n-1}(m)(\mbox{ mod } m^{4l})$, such that
$$g=x_{5n-4}\cdot \cdot\cdot x_{2n-1}\gamma_2 x_{2n-2}\cdot \cdot\cdot x_1 \gamma_1.$$
\epr

In the next proposition, we denote by $\tilde{A}$ the projection of a subset $A\subset \Ga_n(m,l)$ modulo $m^{4l}$.
\bpr \label{pr_kaz_G_n_l} Let $n\geq 3$ and $m\in \mathbb{N}$. then there exists $\mu=\mu(n,m) >0$ such that for any $l\in \mathbb{N}$ we have
$$\kappa(G_n(m,l),\tilde{\beta}_n(m,l))>\mu.$$
\epr
\begin{proof}
\noindent From Proposition \ref{ga_mod_bg} we get that for any element $g\in G_n(m,l)$ there exists $x_1,...,x_{5n-4} \in X_n(m^l)(\mbox{ mod } m^{4l})$ and $\gamma_1,\gamma_2\in \Ga_{n-1}(m) (\mbox{ mod } m^{4l})$, such that
$$g=x_{5n-4}\cdot \cdot\cdot x_{2n-1}\gamma_2 x_{2n-2}\cdot \cdot\cdot x_1 \gamma_1.$$

\noindent In the case $n=3$, by Selberg celebrated result \cite{Se}, there exists $\ep'>0$ (which depends on $\tilde{\Sigma}_2(m) \subset \tilde{\beta}_3(m,l) $) such that $\kappa(\Ga_2(m)/\Ga_2(m^{4l}),\tilde{\Sigma}_2(m)) >\ep'$.\\
\noindent In case $n >3$, from Theorem \ref{thm_kaz_el} (2) there exists $\ep'>0$ such that $\kappa(\Ga_{n-1},\Sigma_{n-1})>\ep'$.
Hence  $\kappa(\Ga_{n-1}(m)/\Ga_{n-1}(m^{4l}),\tilde{\Sigma}_{n-1})>\ep'$.

 Let $(\pi,\mathcal{H})$ be a unitary representation of $G_n(m,l)$,
and assume that there exists a unit vector $v$ which is $(\beta_n(m,l),\mu)$-invariant where $0<\mu<\ep'$ will be determined latter.
Now restrict the representation $\pi$ to the subgroup $\Ga_{n-1}(m)/\Ga_{n-1}(m^{4l})$ and by Lemma \ref{lem_kaz_gp} we get
 $$\|\pi(g)v-v\|\leq 2\frac{\mu}{\ep'},$$  for every $g \in \Ga_{n-1}(m)/\Ga_{n-1}(m^{4l})$.
In addition, by Proposition \ref{pro_kaz_ratio_ga} for any $h\in X_n'(m^l) (\mbox{ mod } m^{4l})$ we have
 $$\|\pi(h)v-v\|\leq \frac{\mu}{\frac{\ep(m)}{2}}=\frac{2\mu}{\ep(m)}.$$

\noindent Let $\lambda=\max\{\frac{2\mu}{\ep'},\frac{2\mu}{\ep(m)}\}$, Hence

$$\|\pi(g)v-v\|\leq \sum \limits_{i=1}^{5n-4}\|\pi(x_i)v-v\|+\sum \limits_{j=1}^2\|\pi(\gamma_j)v-v\|\leq (5n-2)\lambda$$

\noindent So from Lemma \ref{cen_mas} for $\lambda= \frac{1}{5n-2}$, this representation contains a non-zero invariant vector.
Therefore for
$$\mu =\max\{\frac{\ep'}{10n-4},\frac{\ep(m)}{10n-4}\},$$
the result follow.

\end{proof}

\bthm  \label{thm_uni_m_l} Let $n\geq 3$ and $m\in \mathbb{N}$. Then there exist $\ep=\ep(n,m)>0$ such that for any $l\in \mathbb{N}$ we have
$$\kappa(\Ga_n(m,l),\beta_n(m,l))> \ep.$$
\ethm

\begin{proof}

Let $(\pi,\mathcal{H})$ be a unitary representation of $\Ga_n(m,l)$.
Assume that there exists a unit vector $0\neq u \in \mathcal{H}$ which is $(\beta_n(m,l),\ep)$-invariant, where $\ep \leq \frac{\delta_n}{2}$ will be
determined latter.
Then from Proposition \ref{prop_inv_el_m_l}, there exist a non-zero invariant vector  in $\mathcal{H}$ under the action
of $EL_n(m^{2l})$ and in particular under the action of the normal subgroup $\Ga_n(m^{4l})\leq EL_n(R_{m^{2l}})$ (see \cite{Men2,Tit}).
Decompose $\mathcal{H}=\mathcal{H}_0+\mathcal{H}_1$ where $\mathcal{H}_0$ is the invariant subspace under the action of $\Ga_n(m^{4l})$ and $\mathcal{H}_1$ is the orthogonal complement.
Accordingly we have $\pi=\pi_0+\pi_1$ and we write $v=v_0+v_1$.

If $\|v_1\| \geq \frac{1}{2}$, then from Proposition \ref{prop_inv_el_m_l} we obtain that
$$\|\pi(s)v-v\|=\|\pi(s)v_1-v_1\|>\delta_n\|v_1\|\geq \frac{\delta_n}{2},$$
\noindent for some $s\in \beta_n(m,l)$ and we are done with $\ep=\frac{\delta_n}{2}$.

So we assume that $\|v_0\| \geq \frac{1}{2}$. Since $\Ga_n(m^{4l})$ is a normal subgroup, therefore $(\pi_0,\mathcal{H}_0)$ is a unitary representation of $\Gamma_n(m,l)$ which give rise to a unitary representation of the quotient $\Ga_n(m,l)/\Ga_n(m^{4l})=G_n(m,l)$.
From Proposition \ref{pr_kaz_G_n_l} there is $\mu>0$ (independent of $l$) such that
$$\kappa(G_n(m,l),\tilde{\Sigma}_n(m))> \mu.$$

\noindent Therefore, $\ep=\min\{\frac{\mu}{2},\frac{\delta_n}{2}\}$ is a Kazhdan constant of $\Ga_n(m,l)$ with respect to $\beta_n(m,l)$.

\end{proof}

Now we are ready to complete the proof. Let $\Ga$ be a finite index subgroup of $SL_n(\mathbb{Z})$ where $n\geq 3$. From the positive solution for the congruence subgroup problem \cite{BLS,Men1}, there exists $m\in \mathbb{N}$ such that $\Ga_n(m)$ is a (finite index) subgroup of $\Ga$.
For any $l\in \mathbb{N}$ the group $\Gamma(m,m^l)$ is a finite index subgroup of $\Ga$. By Theorem \ref{thm_uni_m_l} we get that the family $\{\Gamma_n(m,m^l)\}_{l=1}^\infty$ has uniform Kazhdan constants with respect to the generating sets $\{\beta_n(m,l)\}_{l=1}^\infty$.
\qed

\section {Acknowledgments.} It is my pleasure to thank Martin Kassabov for identifying an error in an early draft of this paper and whose numerous remarks helped me improve the paper. I am grateful to Dave Morris for letting me publish his proof of the bounded generation property for $EL_n(R_m)$. The insightful discussions with Anderi Jaikin-Zapirain, Alex Lubotzky, Chen Meiri and Yehuda Shalom were valuable to this work.
This paper was initiated while the author was a post-doc at the Hebrew university of Jerusalem. I am grateful for their support. I also would like to thank the Arthur and Rochelle Belfer Institute of Mathematics and Computer Science at the Weizmann Institute, the Computer Science Division at the Open university of Israel and the IH\'{E}S at Bures-sur-Yvette for their support. I acknowledge the support of the ISF grant no. 221/07.

\noindent Department of Mathematics, Weizmann Institute of Science\\
 Rehovot 76100, Israel\\
and\\
Computer Science Division, The Open University of Israel\\
\textit{E-mail address: }uzy.hadad@gmail.com

 \end{document}